\documentclass[11pt]{amsart}
\usepackage{geometry,verbatim}                
\usepackage{graphicx}
\usepackage{amssymb}
\usepackage{epstopdf}
\usepackage{lscape}
\input mssymb.tex
\newtheorem{example}{Example}
\newtheorem{definition}{Definition}
\def\R{\Bbb R}
\def\C{\Bbb C}

\def\spp{\vspace{5pt}\noindent}
\newtheorem{proposition}{Proposition}

\title{Low Energy Clamped Planar Elastica}
\author{Lyle Noakes}
\date{August 17, 2016}                                           
\address{Lyle.Noakes@uwa.edu.au (School of Mathematics \& Statistics, The University of Western Australia, 35 Stirling Highway, Crawley, WA 6009, AUSTRALIA)}
\begin{document}
\maketitle
\section{Continuous Elastica}
The Euclidean plane $E^2$ be the Euclidean plane,  equipped with the Euclidean inner product $\langle ~,~\rangle$, is identified with the complex plane $\C$ in the standard way, with $(0,1)$ corresponding to 
${\bf i}$. The unit circle in $E^2$ is denoted by $S^1$, and we suppose that real numbers $a<b$ are given. 
Given $x_a,x_b\in E^2$, together with $v_a,v_b\in S^{1}\subset E^{2}$, suppose that there exists a $C^\infty$ unit-speed curve $x:[a,b]\rightarrow E^{q+1}$ satisfying 
\begin{equation}\label{bcons}x(a)=x_a,~ \dot x(a)=v_a,~ x(b)=x_b,~\dot x(b)=v_b.\end{equation}
Such a curve has length $L:=b-a$, and is said to be {\em feasible} for $a,b,x_a,x_b,v_a,v_b$. An Euler-Bernoulli (fixed length) {\em elastica} is defined to be a critical point of the elastic energy functional
$$E(x)~:=~\frac{1}{2}\int _a^b\Vert \ddot x(t)\Vert ^2~dt$$
as $x$ varies over feasible curves, where $\Vert ~\Vert$ is the Euclidean norm. The beautiful review \cite{levien} of the study of elastica, from James Bernoulli and Leonhard Euler through to 2008, contains many interesting references including \cite{dan}.  As explained in \cite{levien}, the difficulty of obtaining numerical solutions for elastica satisfying prescribed conditions was influential in the development of the modern theory of splines. The present paper attempts a small additional contribution, by way of simplicity and speed, to the well studied area 
of numerical methods for elastica \cite{edwards}, \cite{bruckstein}.

\spp
From the Pontryagin Maximum Principle \cite{ponty}, a feasible $x$ is an elastica when, for some $C^\infty$ function 
$\mu :[a,b]\rightarrow \R$, 
\begin{equation}\label{eleqctselas}
\ddot v(t)+\mu (t)v(t)~=~C
\end{equation}
where $v(t):=\dot x(t)$, and $C\in E^{q+1}$ is constant. Taking inner products with $v(t)\in S^1$, we see that 
$$\mu (t)~=~\langle C,v(t)\rangle -\langle \ddot v,v\rangle~=~\langle C,v(t)\rangle +\kappa (t)^2$$
where {\em curvature} is defined by $\kappa (t):=\det \left[ v(t) \vert \dot v(t)\right] =\langle v(t),{\bf i}\dot v(t)\rangle $.
\begin{definition} A {\em lifting} of $v:[a,b]\rightarrow S^1$ is 
 a $C^\infty$ function $\theta :[a,b]\rightarrow \R$  
satisfying $(\cos \theta (t),\sin \theta (t))=v(t)$ for all $t\in [a,b]$.\qed
\end{definition}
\noindent  
For any lifting $\theta $ of $v$, $\kappa =\dot \theta =\pm \Vert \dot v\Vert$. Taking inner products of (\ref{eleqctselas}) with $\dot v(t)$, 
$$\langle \ddot v(t),\dot v(t)\rangle ~=~\langle C,\dot v(t)\rangle~\Longrightarrow~
\langle C,v(t)\rangle ~=~\frac{\kappa (t)^2-c}{2}~\Longrightarrow~\mu (t)~=~\frac{3\kappa (t)^2-c}{2},$$
where $c\in \R$ is constant. Differentiating 
 (\ref{eleqctselas}) and taking inner products with $\dot v(t)$,
$$0~=~\langle v^{(3)},\dot v\rangle +\mu \kappa ^2~=~\kappa \ddot \kappa +\dot \kappa ^2-\langle \ddot v,\ddot v\rangle +\mu \kappa ^2~=~\frac{\kappa (t)}{2}(2\ddot \kappa (t)+\kappa ^3(t)-c\kappa (t)  )$$
because $\ddot v=-\kappa ^2v+\dot \kappa {\bf i}v $. Therefore (or, alternatively, following the derivation in \cite{singer}), 
\begin{equation}\label{eleq}2\ddot \kappa (t)~=~c \kappa (t)-\kappa (t)^3.\end{equation}
Excluding the trivial cases where $\kappa $ is constant, namely $x$ is either a circular arc or a line segment, the solutions of (\ref{eleq}) are 
\begin{equation}
\kappa (t)^2~=~\kappa _0^2(1-\frac{p^2}{w^2}{\rm sn} ^2\left( \frac{\kappa _0}{2w}(t-t_0),p\right) )
\end{equation}
where ${\rm sn}$ denotes the {\em elliptic sine}, $w$ is either $p$ or $1$, and $c$ is related to the parameters $\kappa _0,p$ by 
$$2c ~=~\frac{\kappa _0^2}{w^2}(3w^2-p^2-1).$$
For $w=p=1$ the elastica is called {\em borderline}. Otherwise, according as $w=p$ or $w=1$, it is said to be {\em wavelike} or {\em orbitlike}.  
\section{Nontrivial Cases}
If $x$ is an elastica then, for any Euclidean transformation $A$ of $E^2$, so is $t\mapsto Ax(t)$ . 
So suppose, without loss,  that $x_a=(0,0)$ and $v_a=(1,0)$.
\subsection{Wavelike} For a wavelike elastica, $\kappa $ oscillates periodically between $\pm \kappa _0$, according to 
\begin{equation}\label{wl}\kappa (t)~=~\kappa _0{\rm cn}\left(  \frac{\kappa _0}{2p}(t-t_0),p\right) ,\end{equation}
where the {\em elliptic cosine} ${\rm cn}$ is given by ${\rm cn}(u,p)=\cos \phi $ where $\phi $ is the {\em Jacobi amplitude} ${\rm am}(u,p)$, namely 
$$u~=~F(\phi ,p)~:=~\int _0^\phi \frac{1}{\sqrt{1-p\sin ^2\psi}}~d\psi .$$

\subsection{Orbitlike} For an orbitlike elastica we have
\begin{equation}\label{ol}\kappa (t)~=~\kappa _0{\rm dn} \left(  \frac{\kappa _0}{2}(t-t_0),p \right) ,\end{equation} 
where ${\rm dn}(u,p):=\sqrt{1-p\sin ^2\phi}$, with $\phi ={\rm am}(u,p)$ as before. Integrating (\ref{ol}),  
$$\displaystyle{\theta (t)=2\left( {\rm am} (\frac{\kappa _0(t-t_0)}{2},p) -{\rm am} (\frac{\kappa _0(a-t_0)}{2},p) \right) },$$
where $\theta :[a,b]\rightarrow \R$ is the lifting of $v$ with $\theta (a)=0$.

\begin{example}\label{ex1} Taking $a=0$, $b=10$, $\kappa _0=1$, $t_0=1/2$ and $p=2$, we find that $\theta (t)=2{\rm am} (t/2-1/4,2)+0.489774$. The corresponding elastica $x:[0,10]\rightarrow E^2$, shown in Figure \ref{fig1}, 
is found by numerically solving $\dot x(t)=(\cos \theta (t),\sin \theta (t))$. Mathematica's NDSolve takes 0.586 seconds on a 1.7GHz Intel Core i5 Mac with 4GB RAM.
\begin{figure}[h]\label{fig1}
  \begin{center}
    \includegraphics[width=2.5in]{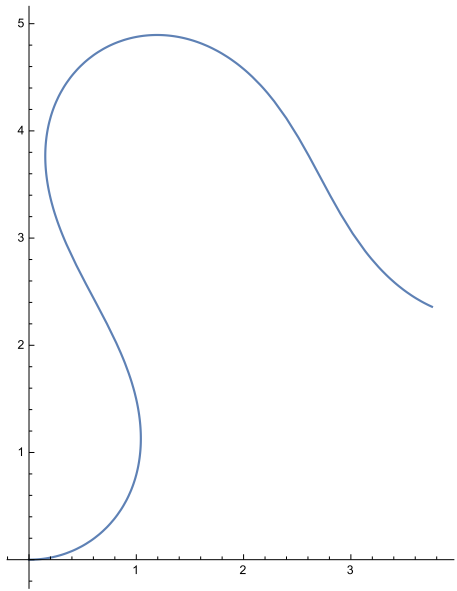}
    \caption{An Orbitlike Elastica (Example \ref{ex1}).}
  \end{center}
  \end{figure}
  \qed
\end{example}
\subsection{Borderline} For a borderline elastica, $\displaystyle{\kappa (t)=\kappa _0{\rm sech} (\frac{\kappa _0}{2}(t-t_0))}$ is nonperiodic. 
\section{Boundary Conditions}\label{bc}
Given $a,b,x_a,v_a$, the elastica $x$ is determined by its curvature $\kappa :[a,b]\rightarrow \R$.  
So the parameters $\kappa _0,p,t_0$ must be chosen 
to satisfy $x(b)=x_b$ and $\dot x(b)=v_b$. This can be done, separately for each case, by numerically solving 
a system of three nonlinear equations, usually in elliptic functions. This is time-consuming and 
there are generally multiple solutions\footnote{  
This contrasts with the analogous problem where $x$ is not necessarily unit-speed and the solution is a unique cubic polynomial.} to the boundary value problem. With applications and extensions in mind, we make it a point to search for uncomplicated clamped elastica, especially {\em minimisers} of $E$.

\spp
Background on finding elastic splines\footnote{A different problem, where $L$ is not considered in advance, is studied by Brunnett and Wendt \cite{brunnett}.} can be found in 
\S16 of \cite{levien}, with more details in \cite{edwards}. More recently, a contribution by Bruckstein, Holta, Netravalia and Arun\cite{bruckstein} solves boundary value and interior value problems of this sort, and in much greater generality than considered in the present paper. Their method, which we call the {\em standard discretisation}, proceeds by 
optimising a discrete analogue of energy for 
piecewise-linear curves satisfying the given constraints. One of the great advantages of standard discretisation is 
ease of implementation, but the method can be time-consuming and may easily 
result in clamped elastica of unnecessarily high energies.

\begin{example}\label{ex2} As in Example \ref{ex1}, take $a=0$, $b=10$, $x_a=(0,0)$ and $v_a=(1,0)$. 
Then (cheating a little) take $x_b=x(b)=(3.75605, 2.35942)$ and $v_b=\dot x(b)=(0.911711, -0.410832)$ where $x$ is the elastica found in Example \ref{ex1}. Subdividing $[0,10]$ into $100$ subintervals, the vertices of the corresponding discrete elastica are shown (red) in Figure \ref{fig2}, together with the original $x$ (blue). 
Standard discretisation takes 597 seconds to find the discrete elastica, which has discrete energy $0.540$ compared with $0.043$ for $101$ equally spaced points along $x$. Although the standard discretisation is straightforward, some appreciable effort is needed to compute it, and the discrete energy is much too high. 
\begin{figure}[h]\label{fig2}
  \begin{center}
    \includegraphics[width=2.5in]{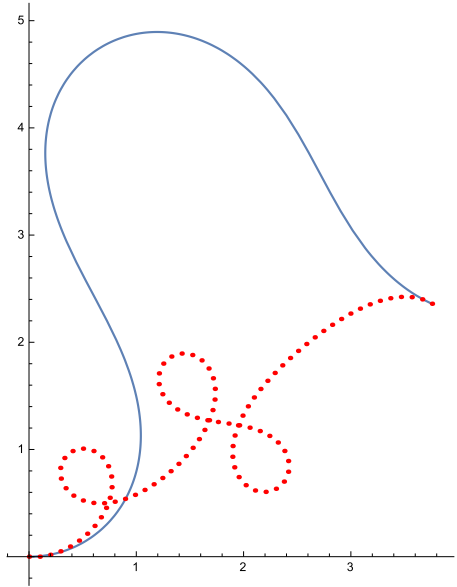}
    \caption{Continuous (blue) and standard discretisation (red)  equal-length elastica, with the same length and end conditions (Examples \ref{ex2}, \ref{ex4}).}
  \end{center}
  \end{figure}
\qed
\end{example} 
Standard discretisation seems a big improvement over classical methods, but cases such as Example \ref{ex2} suggest that it may be   
insufficiently robust for routine use by nonexperts interested in uncomplicated curves. 
%
%
%

\spp
Our method for estimating clamped elastica is designed to search for elastica $x$ of small energy. Indeed, in the first step, we make the stronger assumption that all derivatives of $v$ are moderate in size. 
Then the assumption is used to find a first estimate $\tilde \theta$ of a lifting $\theta$ of the unknown elastica $x$. 
In the second step, $\tilde \theta$ is taken as a starting point for a numerical optimisation of an approximate energy. Because $\tilde \theta$ is already uncomplicated, and reducing energy 
should not make things worse, we are more likely to achieve a global minimum of $E$, 
rather than just a local minimum. 

\spp
Computational speed is addressed in Step 2, where Simpson's Rule improves 
approximations to $E$, and estimates of $\theta$ are modelled as polynomial splines. 
In effect, smoothness of elastica is used to reduce   
the need for a large search space.  Numerical optimisation then proceeds quickly. 

\section{Step 1}\label{step1}
Minimising $E(x)$ means minimising the $L^2$ norm of $\dot v$. 
In this first step we aim for an uncomplicated initial curve $\tilde x$ where the $L^2$ norms of derivatives of all orders are not too large, having regard to $L=b-a$ and approximately satisfying the prescribed boundary conditions. 
Rather than construct $\tilde x$ explicitly, we estimate $\tilde v:=\dot {\tilde x}$, then a lifting $\tilde \theta :[a,b]\rightarrow \R$ of $\tilde v$. Then $\tilde \theta $ is used in Step 2, to start the optimisation of an approximation to $E$. 
\spp
\begin{definition} Writing $h:=(b-a)/4$, a quantity $f$ calculated from $x$ is said to be $O(h^n)$ when, 
for some constant $K$ depending only on some class\footnote{Here $a$ and $b$ may vary, depending 
on the choice of $x$ 
within the class. } to which $x$ belongs (not specifically on $x$), 
the magnitude of $f$ is bounded above by $Kh^n$. So a vector valued function $f$ calculated from $x$ and defined on $[a,b]$ is $O(h^n)$ when, for some constant 
$K$ independent of $x$ and all $t\in [a,b]$, we have $\Vert f(t)\Vert \leq Kh^n$. \qed
\end{definition}

\spp
Assuming  $v^{(m)}=O(1)$ for $1\leq m\leq 4$, write  $v_k:=v(a+kh)$ for $0\leq k\leq 4$ so that, approximating $v^{(4)}$ by central differences,     
$$\displaystyle{\frac{v_4-4v_3+6v_2-4v_1+v_0}{h^4}~=~v^{(4)}(a+2h)+O(h^2)~=~O(1)\Longrightarrow}$$ 
\begin{equation}\label{v4est}
-2(v_1+v_3)+3v_2~=~-\frac{1}{2}(v_a+v_b)+O(h^4). 
\end{equation}
From $\int _a^bv(t)~dt~=~x_b-x_a$ we find, using\footnote{Boole's Rule might be used instead, or (with additional complexity) Gauss-Lobatto quadrature, but there is a lot to be said for simplicity.} the Composite Simpson's Rule,  
\begin{equation}\label{quadeq}
2(v_1+v_3)+v_2~=~\frac{3}{2h}(x_b-x_a)-\frac{1}{2}(v_a+v_b)+O(h^{4}).
\end{equation}
Eliminating $v_1+v_3$ between (\ref{v4est}), (\ref{quadeq}), we find that $v_2=\tilde w_2+O(h^4)$ where 
\begin{equation}\label{w2til} \tilde w_2~:=~\frac{3}{8h}(x_b-x_a)-\frac{1}{4}(v_a+v_b).\end{equation}
Therefore, and because $\Vert v_2\Vert =1$, we have $v_2=\tilde v_2+O(h^4)$, where $\tilde v_2:=\tilde w_2/\Vert \tilde w_2\Vert$. 

\spp
Substituting for $v_2$ in (\ref{quadeq}), we find that $2w_{1,3}:=v_1+v_3=2\tilde w_{1,3}+O(h^4)$ where
\begin{equation}\label{quadeqs}
\tilde w_{1,3}~:=~\frac{3}{8h}(x_b-x_a)-\frac{1}{8}(v_a+v_b)-\frac{1}{4}\tilde v_2.
\end{equation}
By Taylor's Formula\footnote{Other kinds of estimates can also be made, but this has the virtue of simplicity.},   
\begin{equation}\label{tf}v_3-v_1~=~\frac{1}{2}(v_b-v_a)+O(h^3),\end{equation}
so that $v_1=\tilde w_1+O(h^3)$ and $v_3=\tilde w_3+O(h^3)$ where 
\begin{equation}\label{weqs}
\tilde w_1~:=~\tilde w_{1,3}-\frac{1}{4}(v_b-v_a) \quad \hbox{and}\quad \tilde w_3~:=~\tilde w_{1,3}+\frac{1}{4}(v_b-v_a).
\end{equation}
Then, because $v_1$ and $v_3$ are unit vectors, $v_1=\tilde v_1+O(h^3)$ and $v_3=\tilde v_3+O(h^3)$ 
where $\tilde v_j:=\tilde w_j/\Vert \tilde w_j\Vert$. So we have estimated $v_1$ and $v_3$ to $O(h^3)$ errors, 
and $v_2$ to $O(h^4)$. Summarising so far,
\begin{proposition}\label{first}  Given $a<b\in \R$, $x_a,x_b\in E^2$, $v_a,v_b\in S^1$, define $\tilde v_j:=w_j/\Vert w_j\Vert \in S^1$ for $j=1,2,3$, where $w_2$ is given by formula (\ref{w2til}), and $w_1,w_3$ by (\ref{weqs}). 
Then for suitably small $h$, and assuming that derivatives of $x$ are $O(1)$, we have 
$v_2=\tilde v_2+O(h^4)$ and, for $j=1,3$, $v_j=\tilde v_j+O(h^3)$. \qed
\end{proposition}

\spp
Next the $\tilde v_j$ are used to find a rough estimate $\tilde \theta $ of 
the lifting $\theta :[a,b]\rightarrow \R$ of $v$  where,  
without loss of generality, $v_a=(1,0)$ with $\theta (a)=0$. 
We write $\tilde v_4:=v_b$, $\tilde \theta _0:=\tilde \theta (a)=0$ and require
\begin{equation}\label{thcons}(\cos \tilde \theta (a+jh),\sin \theta (a+jh))~=~\tilde v_j\quad \hbox{for}\quad 1\leq j\leq 4.\end{equation}

\spp
For $1\leq j\leq 4$, this only determines the $\tilde \theta _j:=\tilde \theta (a+jh)$ modulo $2\pi$. To encourage simpler elastica 
the $\tilde \theta _j$ are chosen\footnote{This corresponds to smallest total curvatures over the $[a+(j-1)h,a+jh]$ for $j=1,2,3,4$.}  for $j=1,2,3,4$, as close as possible to $\tilde \theta _{j-1}$ consistent with (\ref{thcons}). Then we interpolate accordingly\footnote{The specific method of interpolation is not especially critical. Much more significantly, our construction of $\tilde \theta $ discourages exotic solutions of the subsequent optimisation problem.}. 
A reasonable choice for $\tilde \theta :[a,b]\rightarrow \R$ is the natural cubic spline\footnote{This corresponds to minimising the $L^2$ norm of $\dot \kappa$, 
at least for the initial guess.} satisfying 
$\tilde \theta (a+jh)=\tilde \theta _j$ for $0\leq j\leq 4$. 

\spp
In short, Step 1 proceeds as follows:
\begin{enumerate}
\item If necessary, translate and rotate the data so that $x_a=(0,0)$ and $v_a=(0,1)$. Set $h:=(b-a)/4,\tilde v_0:=v_a,\tilde v_4:=v_b$. 
\item Define $\tilde w_2$ is given by formula (\ref{w2til}). 
\item Define $\tilde w_{1,3}$ by formula (\ref{quadeqs}), then $\tilde w_1$ and $\tilde w_3$ by formula (\ref{weqs}). 
\item Set $\tilde v_j:=\tilde w_j/\Vert \tilde w_j\Vert$ for j=1,2,3.
\item Set $\tilde \theta _0:=0$ and, for $1\leq j\leq 4$, choose $\tilde \theta _j$ so that $(\cos \tilde \theta _j,\sin \tilde \theta _j)=\tilde v_j$ with the $\vert \tilde \theta _j-\tilde \theta _{j-1}\vert $ as small as possible. 
\item Let $\tilde \theta :[a,b]\rightarrow \R$ be the natural cubic spline satsifying $\tilde \theta (a+jh)=\tilde \theta _j$ for $0\leq j\leq 4$. 
\end{enumerate}

\spp
Although $\tilde \theta $ approximates $\theta $ with at most $O(h^3)$ errors, the actual bounds on the $v^{(m)}$ may be difficult to estimate. So in practice it may be hard to say exactly how good the approximation is. 
Our algorithm is intended for relatively uncomplicated elastica $x$. So it is interesting to compare $\theta $ and $\tilde \theta $ in Example \ref{ex1} which, as seen in Example \ref{ex2}, is a nontrivial case. 
\begin{example}\label{ex3} In Figure \ref{fig3}, the graph of $\tilde \theta$ (yellow) is not a highly accurate 
estimate of the graph of $\theta$ (blue) for Example \ref{ex1}. On the other hand, there do not seem to be any very remarkable differences between the two curves: this is all that is needed to begin the second step of our algorithm. It takes 56 seconds to plot the graph of $\theta$ by integrating the known closed-form solution for $\kappa$, compared with $0.197$ seconds for plotting $\tilde \theta$. So solving the initial value problem from the closed-form solution is already time-consuming. Our algorithm, whose second step is given in Section \ref{nextstep}, solves the much harder boundary-value problem. 
\begin{figure}[h]\label{fig3}
  \begin{center}
    \includegraphics[width=2.5in]{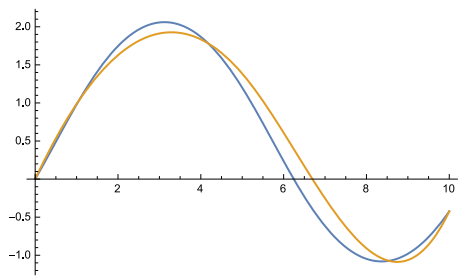}
    \caption{Liftings $\theta$ (blue) and $\tilde \theta $ (yellow) in Example \ref{ex3}.}
  \end{center}
  \end{figure}
  \qed
  \end{example}
  \section{Step 2}\label{nextstep}
Step 1 gives a lifting $\tilde \theta$ of some uncomplicated curve $\tilde x$ 
that approximately satisfies the end conditions.  Next we approximately minimise the elastic energy, and more nearly satisfy the end conditions. This is done as follows.

\begin{enumerate}
\item For an integer $n$ greater than $4$, redefine $h:=(b-a(/n$ and evaluate 
$\tilde \theta$ at $a+jh$ for $0\leq j\leq n$ to give $\phi _{(0)}\in \R ^{n+1}$. 
\item For variable\footnote{corresponding notionally to an unknown lifting $\theta$ of $v$ for the unknown uncomplicated elastica $x$} $\phi \in \R ^{n+1}$, approximate $\int _a^bv(t)dt\in E^2$ by a sum of the form 
~$\displaystyle{S(\phi ):=\sum_{j=0}^nq_j(\cos \phi _j,\sin \phi _j)\in E^2}$\\ 
for suitable constant\footnote{When $n$ is even, the Composite Simpson's Rule may be used,   
namely $q=h(1,4,2,4,2,4,\ldots ,2,4,1)/3$.} $q=(q_0,q_1,\ldots ,q_n)\in \R^{n+1}$. 
\item Starting with $\phi _{(0)}$ as an initial guess, numerically minimise $\sum_{j=1}^n(\phi _j-\phi _{j-1})^2$ 
with $\phi _0=\theta _a$, subject to the trigonometric constraints $S(\phi )=x_b-x_a$.
\item Interpolate the\footnote{we hope, but possibly local minimiser} minimiser $\phi$ by some convenient\footnote{such as the natural cubic spline} $C^2$ curve $\hat \theta :[a,b]\rightarrow \R$. 
\item Setting $\hat v(t):=(\cos \hat \theta (t),\sin \hat \theta (t))$, take $\hat x$ to be a numerical solution\footnote{Mathematica's {\em NDSolve} may be used.} of 
$\dot {\hat x}(t) =\hat v(t)$ with $\hat x(a)=x_a$.   
\end{enumerate}
There is another option, namely to start Step 2 with with $n$ small, then gradually increase $n$, repeating Step 2 with initial estimates of $\theta$ from previous optimisations. 
This more gradual movement from $\tilde x$ to $\hat x$ might occasionally be advantageous, but   
we have had no difficulty with the method as it stands. 
 As explained in \S \ref{bc} and illustrated in \S \ref{exsec}, 
 it seems reasonable to hope that the present method for finding clamped elastica is more robust, gives better results,  and is faster than standard discretisation. 

\section{Comparisons with Standard Discretisation}\label{exsec}
\begin{example}\label{ex4} With boundary conditions from Example \ref{ex2} and $n=20$, the estimate $\hat x:[a,b]\rightarrow E^2$ is almost\footnote{With $n=40$ there is no visible difference at all.} indistinguishable from the original elastica $x$ in Figure \ref{fig2}. It takes $0.066$ seconds to compute $\hat x$, compared with 
$597$ seconds for the standard discretisation (dotted). The discrete energy of $\hat x$ is $0.0425$, compared with 
$0.54$ for the standard discretisation.  
 \qed
\end{example}
\begin{example}\label{ex5} Taking $b=15$ instead of $10$ in Example \ref{ex4}, we have  
$x_b=(4.38081, 6.00329)$, $v_b=(-0.0106571, 0.999943)$, the standard discreisation (red in Figure \ref{fig3}) takes $288.54$ seconds. The standard discretisation has energy $0.0968$ compared with $0.1465$ 
for the original elastica $x$ and, correspondingly, has a somewhat simpler appearance. So, on this occasion, the standard discretisation is preferable\footnote{The aim is not to recover $x$, rather to minimise elastic energy subject to the given length and boundary conditions (in this case read from $x$).} to $x$.  Our method is better still, taking $0.054$ 
seconds to compute $\hat x$. The even less complicated appearance of $\hat x$ (continuous red) in Figure \ref{fig3} is consistent with its lower discrete energy of $0.0292$. \qed
\begin{figure}[h]\label{fig4}
  \begin{center}
    \includegraphics[width=2.5in]{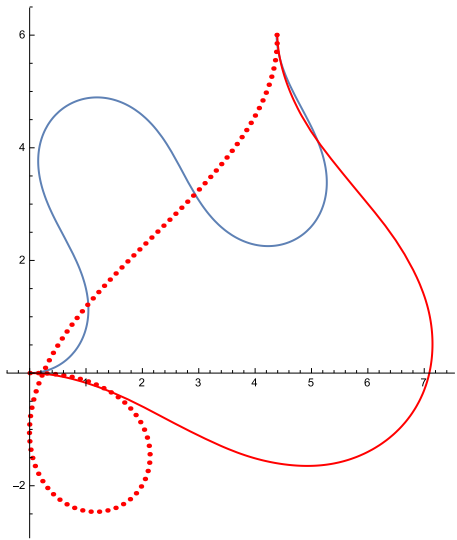}
    \caption{$x$ (blue), standard discretisation (red dotted) and $\hat x$ (continuous red) in Example \ref{ex5}.}
  \end{center}
  \end{figure}
\end{example}
\begin{example}\label{ex6} Increasing $b$ to $20$, $x$ becomes more convoluted with discrete energy $0.351$. Using boundary data from $x$, Mathematica takes 
$778$ seconds to report failure of standard discretisation (nonconvergence, unusable output). So there is some question about robustness of standard discretisation, at least when pushed to this extent. 
On the other hand, our algorithm for finding  $\hat x$ takes $0.122$ seconds. Consistent with its uncomplicated appearance (red in Figure \ref{fig5}), $\hat x$ has discrete energy 
$0.0758$, which compares well with the energy of $x$.  \qed
\begin{figure}[h]\label{fig5}
  \begin{center}
    \includegraphics[width=2.5in]{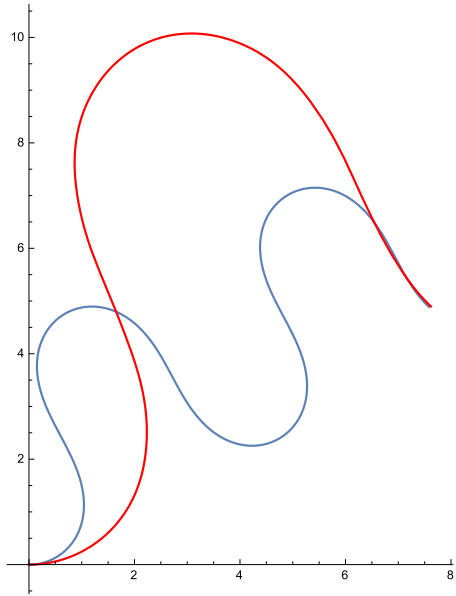}
    \caption{$x$ (blue) and $\hat x$ (red) in Example \ref{ex6}.}
  \end{center}
  \end{figure}
\end{example}

\section{Conclusion}
A method is given for estimating clamped plane elastica. Arguments are made, and evidence is provided by way of illustrative examples,  suggesting that the new method is quicker and more robust than standard discretisation, and more likely to give elastica of low energy. Just as for standard discretisation, no use is made of the known solutions for elastica in terms of elliptic functions.  
An extension to calculating general elastic splines is kept for a future paper.

\end{document}